# Spatial modelling for mixed-state observations

## Cécile Hardouin

*SAMOS-MATISSE-Centre d'Économie de la Sorbonne*
*Université de Paris 1*
*90, rue de Tolbiac*
*75634 Paris Cedex 13*
*e-mail:* hardouin@univ-paris1.fr

## Jian-Feng Yao

*IRMAR*
*Université de Rennes 1*
*Campus de Beaulieu*
*35042 Rennes Cedex, France*
*e-mail:* jian-feng.yao@univ-rennes1.fr

**Abstract:** In several application fields like daily pluviometry data modelling, or motion analysis from image sequences, observations contain two components of different nature. A first part is made with discrete values accounting for some symbolic information and a second part records a continuous (real-valued) measurement. We call such type of observations "mixed-state observations".

This paper introduces spatial models suited for the analysis of these kinds of data. We consider multi-parameter auto-models whose local conditional distributions belong to a mixed state exponential family. Specific examples with exponential distributions are detailed, and we present some experimental results for modelling motion measurements from video sequences.

Received January 2008.

## 1. Introduction

In many applications, it is frequent to get observations with two components of a different nature: the first component is made up of discrete values and the second component records a continuous measurement. For example, pluviometry time series at a given site records many zeros for dry days, followed by positive and continuous records for wet periods [2, 1]. Similar phenomena also occur in speech recordings where interchanges are permanent between absences and presences of the signal. Other examples arise in the motion analysis problem from image sequences [5], or in epidemiological data analysis where the disease at given locations can be absent or spreads out. We call such type of measurements *mixed-state observations*. It then raises the question to find accurate models for these types of data.

To deal with data of mixed nature, most of the existing approaches rely on an hierarchic approach. One introduces a hidden variable to distinguish discrete observations from continuous ones. Or equivalently, the discrete values are interpreted as resulting from some unobserved censoring variable [2]. Specifically, a Bayesian approach is used for statistical inference.

Our approach is different. We propose a direct modelling by considering random variables which can take discrete values as well as continuous ones. Although the idea seems absolutely natural, we are not aware of any statistical models relying on such a direct approach for mixed-state data.

The main motivation of the paper is a search for spatial models for observations $\{X_s\}$ such that each $X_s$ is a mixed-state random variable. In the spatial context, the discrete components could not be simply neglected, because these symbolic values as well as their spatial correlations convey important point-wise and contextual information. To this end, we introduce a new class of auto-models for such mixed-state data. Their construction proposed in this paper relies on an adaptation to the present context of a general class of Markov random fields models, namely multi-parameter auto-models, that we recently introduced in [9]. Roughly speaking, an auto-model, as introduced originally in [4], is a Markov field on a finite set of sites, for which the interactions between sites are pairwise only, and each local conditional distributions belongs to some exponential family. The multi-parameter auto-models of [9] extend the classical one-parameter auto-models of [4] and several known spatial models previously proposed in [6, 10, 11].

The plan of the paper is the following. We first present mixed state random variables in a simple context where the observation is made up with 0 and values in $(0, \infty)$. The distribution of this mixed state random variable has two main features; it reflects the dual character of the observation, and the distribution belongs to an exponential family. In §3, we give the general definition for mixed-state variables. We recall in §4 results on general multi-parameters auto-models of [9] which constitute the starting blocks of our construction of auto-models for mixed-state observations that we present in §5. We wet out in §6 a detailed study of mixed-state auto-models where neighbouring sites are spatially cooperative. This property contrasts with many classical auto-models introduced in [4] which lead to a spatial competitive behaviour, which is clearly inadequate in many practical situations. We conclude the paper by an analysis of motion measurements from video sequences, using a mixed positive Gaussian auto-model.

## 2. Simple random variables with mixed states

Before defining general mixed state variables, let us begin with the simplest situation where the state space is $E = \{0\} \cup (0, \infty)$. Of course $E = [0, \infty)$, but the split formula has the merits to insist on the null value which plays a particular role in the construction. A *mixed-state random variable* $X$ on $E$ is defined as follows: with a probability $\gamma \in [0,1]$ we set $X = 0$, and with probability $1-\gamma$, $X$ follows a continuous distribution on $(0, \infty)$ with a probability density function $g$.

Formally, we equip $E$ with its Borel field $\mathscr{E}$ and we introduce a reference measure of mixture type

$$m(dx) = \delta_0(dx) + \lambda(dx), \tag{2.1}$$

where $\delta_0$ is the Dirac measure at 0 and $\lambda$ the Lebesgue measure. Throughout the paper, we denote by $1_A$ the indicator function of a set $A$. For the particular case of $\{0\}$, we use a simpler notation by setting $\delta(x) = 1_{\{0\}}(x)$ and $\delta^*(x) = 1 - \delta(x)$. The above mixed-state variable $X$ then has a probability density function with respect to $m$ given by

$$f(x) = \gamma\delta(x) + (1-\gamma)\delta^*(x)g(x), \quad x \in E. \tag{2.2}$$

Clearly, such mixed-state random variables (or distributions) can provide accurate modelling for the marginal empirical distributions discussed in §1.

For the upcoming construction of spatial models, we are interested in mixed-state random variables of a particular type, namely their continuous component $g$ belongs to a $\ell$−dimensional exponential family

$$g(x) = g_\xi(x) = H(\xi)L(x)\exp\langle \xi, T(x) \rangle, \quad \xi \in \mathbb{R}^\ell, \quad T(x) \in \mathbb{R}^\ell, \tag{2.3}$$

for some sufficient statistics $T$ and measurable positive functions $H$ and $L$ ($\langle,\rangle$ denoting the scalar product in $\mathbb{R}^\ell$). Interestingly enough, the mixed-state distribution can also be put in the form of an exponential family. Indeed,

$$\begin{aligned} f(x) = f_\theta(x) &= \gamma\delta(x) + (1-\gamma)\delta^*(x)g_\xi(x) \\ &= \exp\left[\delta^*(x)\ln\frac{(1-\gamma)H(\xi)}{\gamma} + \langle \xi, T(x)\delta^*(x)\rangle + \log\gamma + \delta^*(x)\log L(x)\right] \\ &= H'(\theta)L'(x)\exp\langle \theta, B(x)\rangle, \end{aligned} \tag{2.4}$$

with $H'(\theta) = \gamma$, $L'(x) = \exp\{\delta^*(x)\log L(x)\}$, and the natural parameter and the sufficient statistics defined by

$$\theta = \begin{pmatrix}\theta_1 \\ \theta_2\end{pmatrix} = \begin{pmatrix}\log\frac{(1-\gamma)H(\xi)}{\gamma} \\ \xi\end{pmatrix}, \quad B(x) = \begin{pmatrix}\delta^*(x) \\ T(x)\delta^*(x)\end{pmatrix}, \quad x \in E. \tag{2.5}$$

Note that with the standard convention $0\log 0 = 0$, these formula are still valid in the extreme situations $\gamma \in \{0,1\}$ which correspond to a purely continuous and a purely discrete distribution, respectively. Therefore, the mixed-state distribution $f_\theta$ belongs to an exponential family, of dimension $\ell + 1$. Moreover the original parameters $\xi$ and $\gamma$ can be recovered from $\theta$ by

$$\xi = \theta_2, \quad \gamma = \frac{H(\theta_2)}{H(\theta_2) + e^{\theta_1}}.$$

Let us consider some examples.

*Example* 1.  Mixed-state Exponential distribution: this simple distribution is obtained with $g_\lambda(x) = \lambda e^{-\lambda x}$ where $\lambda > 0$. Then $\xi = H(\xi) = \lambda$ and $T(x) = -x$. The parametric dimension of the resulting mixed-state distribution equals two.

*Example* 2. Mixed-state Gamma distribution: this situation generalises Example 1 by substituting a Gamma distribution $\Gamma(a,b)$, $a, b > 0$, for the exponential distribution. Here we have $\xi = (b, a-1)$, $H(\xi) = \Gamma(a)^{-1} b^a$ and $T(x) = (-x, \ln x)$. The resulting mixed-state distribution belongs to an exponential family of dimension three.

*Example* 3.  Positive mixed-state Gaussian distribution: here the continuous component of $X$ is the distribution of the modulus of a zero-mean Gaussian distribution with variance $\sigma^2$. We have $\theta = (\ln \frac{2(1-\gamma)}{\gamma \sigma \sqrt{2\pi}}, \frac{1}{2\sigma^2})^\mathrm{T}$ and $B(x) = (\delta^*(x), -x^2)^\mathrm{T}$.

## 3. General random variables with mixed states

To cover situations involving several atomic values, the previous simple model need to be extended. Let $F = \{e_1, \ldots, e_M\}$ be a finite set of $M$ elements and $G$ a Borel subset of an Euclidean space $\mathbb{R}^p$. Let $\mathbf{q} = (q_1, \ldots, q_M)$ be a probability distribution on $F$ and $g$ a probability density function on $G$ (with respect to the Lebesgue measure).

We define a *general mixed-state random variable* $X$ as follows:

- with a probability $\gamma \in [0,1]$, $X$ takes values in $F$ with distribution $\mathbf{q}$;
- with probability $1 - \gamma$, $X$ takes values in $G$ according to the density function $g$.

Although the nature of the discrete state space $F$ could be arbitrary (possibly qualitative), we are going to embed $F$ in $\mathbb{R}^p$ to ease the development of a likelihood-based estimation theory. In other words, we set the state space of $X$ to be
$$E = \{e_1, \ldots, e_M\} \cup G, \quad e_i \in \mathbb{R}^p \setminus G, \quad G \subset \mathbb{R}^p.$$

Therefore, we can supply $E$ with its Borel field and a reference measure of mixture type
$$m(dx) = \sum_{i=1}^M \delta_{e_i}(dx) + \lambda(dx) \tag{3.1}$$

where $\delta_e$ denotes the Dirac measure at a point $e \in \mathbb{R}^p$ and $\lambda$ the Lebesgue measure on $E$. Then the mixed-state variable $X$ has the the following density function with respect to $m$:
$$f(x) = \gamma 1_F(x) \sum_{i=1}^M q_i 1_{\{e_i\}}(x) + (1-\gamma) 1_G(x) g(x), \quad x \in E. \tag{3.2}$$

### 3.1. *Exponential family case*

We focus now on mixed-state distributions which belong to some exponential family. To avoid trivial situations, we assume that

- the discrete distribution **q** is everywhere positive: $q_i > 0$ $(i = 1, \ldots, M)$.
- the density $g$ of the continuous component belongs to a $\ell-$dimensional exponential family as in (2.3).

We first write the discrete distribution **q** in an exponential family form, through the logistic transformation:

$$k_i = \log \frac{q_i}{q_M}, \quad i = 1, \ldots, M.$$

We notice that by definition $k_M = 0$. We have,

$$\sum_{i=1}^{M} q_i 1_{\{e_i\}}(x) = \exp \sum_{i=1}^{M} 1_{\{e_i\}}(x)(k_i + \log q_M), \quad x \in \{e_1, \ldots, e_M\}.$$

Combining this writing with (2.3) and (3.2), we get

$$\begin{aligned}
f(x) &= f_\theta(x) \\
&= \exp\left[1_F(x) \left\{ \sum_{i=1}^{M} 1_{\{e_i\}}(x)(k_i + \log q_M) + \log \gamma \right\} \right. \\
&\quad \left. + 1_G(x) \left\{ \log[(1-\gamma)H(\xi)L(x)] + \langle \xi, T(x) \rangle \right\} \right] \\
&= \exp\left[ 1_G(x) \log \frac{(1-\gamma)H(\xi)}{\gamma q_M} + \langle \xi, T(x) 1_G(x) \rangle \right. \\
&\quad \left. + \sum_{i=1}^{M-1} k_i 1_{\{e_i\}}(x) + \log(\gamma q_M) + 1_G(x) \log L(x) \right] \\
&= H'(\theta) L'(x) \exp\langle \theta, B(x) \rangle, \quad x \in E, \qquad (3.3)
\end{aligned}$$

where $H'(\theta) = \log(\gamma q_M)$, and $L'(x) = L(x)^{1_G(x)}$. In other words, $f_\theta$ belongs to an exponential family of dimension $\ell + M$ with the natural parameter and the sufficient statistics given by

$$\theta = \begin{pmatrix} \theta_1 \\ \vdots \\ \theta_{M-1} \\ \theta_M \\ \theta_{M+1} \end{pmatrix} = \begin{pmatrix} k_1 \\ \vdots \\ k_{M-1} \\ \log \frac{(1-\gamma)H(\xi)}{\gamma q_M} \\ \xi \end{pmatrix}, \quad B(x) = \begin{pmatrix} 1_{\{e_1\}}(x) \\ \vdots \\ 1_{\{e_{M-1}\}}(x) \\ 1_G(x) \\ T(x) 1_G(x) \end{pmatrix}, \quad x \in E. \quad (3.4)$$

Note that $\theta_{M+1}$ and $T(x)1_G(x)$ are $\ell$-dimensional vectors, and by definition $B(e_M) = 0$ and $L'(e_M) = 1$. Furthermore, the original parameters $\xi$, $\mathbf{q}$ and $\gamma$ are recovered through the formulae

$$\begin{aligned}
k_i &= \theta_i, \quad 1 \leq i < M, \\
\xi &= \theta_{M+1}, \\
q_i &= \frac{e^{k_i}}{e^{k_1} + \cdots + e^{k_M}}, \quad 1 \leq i \leq M, \\
\gamma &= \frac{H(\xi)}{H(\xi) + q_M e^{\theta_1}}.
\end{aligned}$$

It is worth noticing that in the case of $E = \{0\} \cup (0, \infty)$, the formulae (3.3)–(3.4) reduce to equations (2.4)–(2.5) of the previous section.

### 3.2. Example of a mixed-state and censored exponential variable

Let $Z$ be an exponential random variable with parameter $\lambda$, censored at a known location $K > 0$. The probability density function of $Z$ is defined by $g(z) = \lambda e^{-\lambda z} \, 1_{(0,K)}(z) + e^{-\lambda K} \delta_K(z)$.

We define the following mixed-state variable $X$: with probability $\alpha$, $X$ takes the value 0 ; with probability $1 - \alpha$, $X$ has the distribution of $Z$. Therefore, $X$ has masses $\{\alpha, (1-\alpha)e^{-\lambda K}\}$ on the atoms $\{0, K\}$, and a continuous density function $(1-\alpha)\lambda e^{-\lambda x}$ on $(0, K)$. Equivalently, $X$ can be viewed as a general mixed-state variable with the state space $E = \{0, K\} \cup (0, K) = [0, K]$ and the following parameters

$$\begin{aligned}
\gamma &= \alpha + (1-\alpha)e^{-\lambda K}, \\
\mathbf{q} &= \frac{1}{\alpha + (1-\alpha)e^{-\lambda K}}(\alpha, (1-\alpha)e^{-\lambda K}), \\
g_\lambda(x) &= \frac{\lambda}{1 - e^{-\lambda K}} e^{-\lambda x}, \quad x \in (0, K).
\end{aligned}$$

Following (3.3) and (3.4), the distribution of $X$ belongs to an 3-dimensional exponential family with the following natural parameters and sufficient statistics:

$$\theta = \begin{pmatrix} \log \frac{\alpha}{1-\alpha} + \lambda K \\ \log \lambda + \lambda K \\ \lambda \end{pmatrix}, \quad B(x) = \begin{pmatrix} 1_{\{0\}}(x) \\ 1_{(0,K)}(x) \\ -x 1_{(0,K)}(x) \end{pmatrix}, \quad x \in E.$$

## 4. Results on multi-parameter auto-models

The construction of spatial models for mixed-state observations relies on the general theory of multi-parameter auto-models developed in [9]. We quote below its main results which are relevant for the present purpose.

Let us set some notations. Let $S = \{1, \ldots, n\}$ be a finite set of sites equipped with a symmetric graph $\mathcal{G}$ without loops. We denote by $\{i, j\}$ a pair of neighbouring sites (in particular, $i \neq j$). For any subset $A$ of $S$, let $x_A = (x_i, i \in A)$ and $x^A = (x_j, j \in S \setminus A)$. The neighbourhood of a site $i$ is $\partial i = \{j \in S : \langle i, j \rangle \}$. We shall write $x^i = x^{\{i\}}$. The variates $x_i$ take their values in a measurable state space $(E, \mathcal{E}, m)$. Most of the time, $E$ will be a subset of $\mathbb{R}^p$. The configuration space $\Omega = E^S$ is equipped with the $\sigma$-algebra $\mathcal{E}^{\otimes S}$ and the product measure $\nu := m^{\otimes S}$. A random field is specified by a probability distribution $\mu$ on $\Omega$, and we assume the positivity condition, that is $\mu$ has an everywhere positive density $P$ with respect to $\nu$. Consequently we write

$$\mu(dx) = P(x)\nu(dx), \qquad P(x) = Z^{-1} \exp Q(x), \tag{4.1}$$

where $Z$ is a normalisation constant. From the Hammersley-Clifford Theorem, the energy function $Q(x)$ is a sum of potentials $\{G_A; A \in \mathscr{C}\}$ indexed by a set $\mathscr{C}$ of cliques. Let us fix a reference configuration, or "ground states", $\tau = (\tau_i) \in \Omega$ yielding to the potentials normalisation: for any potential $G_A(x_A)$ where $A \subset S$, we have $G_A(x_A) = 0$ if $x_i = \tau_i$ for some $i \in A$. This implies $Q(\tau) = 0$ and $Z^{-1} = P(\tau)$ in (4.1).

The class of multi-parameter auto-models defined in [9] extends the classical one-parameter auto-models of J. Besag in its seminal paper [4]. Their construction rely on the following assumptions.

*Assumption* 1. The dependence between the sites is pairwise only; in other words,

$$Q(x) = \sum_{i \in S} G_i(x_i) + \sum_{\{i,j\}} G_{ij}(x_i, x_j).$$

*Assumption* 2. For an integer $k \geq 1$ and all $i \in S$, the conditional distribution of $X_i$ given $X^i = x^i$ relies in an exponential family

$$\log f_i(x_i|x^i) = \langle \theta_i(x^i), B_i(x_i) \rangle + C_i(x_i) + D_i(x^i), \quad \theta_i(x^i) \in \mathbb{R}^k, \quad B_i(x_i) \in \mathbb{R}^k.$$

*Assumption* 3. For all $i \in S$, $\mathrm{Span}\{B_i(x_i) : x_i \in E\} = \mathbb{R}^k$.

The following result of [9] determines the necessary form of the local natural parameters $\{\theta_i(.)\}$ to ensure the compatibility of the family of conditional distributions.

**Theorem 1.** *(Hardouin and Yao [9]). Assume that Assumptions 1–3 are satisfied with the normalisation $B_i(\tau_i) = C_i(\tau_i) = 0$ in Assumption 2 for all $i \in S$. Then, necessarily, the functions $\theta_i$ take the form*

$$\theta_i(x^i) = \alpha_i + \sum_{j \neq i} \beta_{ij} B_j(x_j), \quad i \in S, \tag{4.2}$$

*where $\{\alpha_i : i \in S\}$ is a family of $k$-dimensional vectors, and $\{\beta_{ij} : i, j \in S, i \neq j\}$ is a family of $k \times k$ matrices $\{\beta_{ij}\}$ satisfying $\beta_{ij}^{\mathrm{T}} = \beta_{ji}$. Moreover, the potentials are given by*

$$G_i(x_i) = \langle \alpha_i, B_i(x_i) \rangle + C_i(x_i), \tag{4.3}$$

$$G_{ij}(x_i, x_j) = B_i^{\mathrm{T}}(x_i) \beta_{ij} B_j(x_j). \tag{4.4}$$

A model satisfying the assumptions of this theorem is called a *multi-parameter auto-model*. For a concrete construction of such a multi-parameter auto-model, one follows a two steps method: first, specify the family of cliques (Assumption 1), and the family of conditional distributions (Assumption 2); secondly, find the admissible set of parameters $\{\alpha_i, \beta_{ij}\}$ which ensures the integrability condition:

$$\int_\Omega e^{Q(x)} \nu(dx) < \infty. \tag{4.5}$$

We refer the reader to [9] and the references therein for more account on this new family of auto-models.

Another important question about the model is that of spatial symmetry. The general formulation given above does not impose any symmetry, and hence it can be useful for modelling random fields on arbitrary or oriented graphs. However, in the case of a spatially symmetric random field, all potentials $G_{ij}(x_i, x_j)$ are necessarily symmetric functions; equivalently, all the matrices $\beta_{ij}$ are symmetric.

## 5. Auto-models for mixed state variables

### *5.1. The construction*

Following the general theory of multi-parameter auto-models quoted above, we now construct auto-models for general mixed-state variables $X = \{X_i, \ i \in S\}$ on a finite set $S$. The state space for each variable $X_i$ is $E = F \cup G$ as defined in §3. We let the configuration space $\Omega = E^S = (F \cup G)^S$ be supplied with the product measure $\nu = m^{\otimes S}$, where $m$ is defined by (3.1).

We assume that Assumptions 1–3 are satisfied, where in Assumption 2, the family of conditional distributions $f_i(x_i \mid \cdot)$ belongs to the family of mixed-state distribution given in (3.3). In other words,

$$f_i(x_i \mid x^i) = H'_i(x^i) L'_i(x_i) \exp\langle \theta_i(x^i), B_i(x_i)\rangle \tag{5.1}$$

with $\theta_i(\cdot) \in \mathbb{R}^{\ell+M}$ and $B_i(\cdot) \in \mathbb{R}^{\ell+M}$ satisfying (3.4).

Note that by definition, $B_i(e_M) = 0$ and $L'(e_M) = 1$. Therefore, the state $e_M$ serves as a reference state for coordinates $X_i$ and the reference configuration becomes $\tau = (e_M, \ldots, e_M)$ for the application of Theorem 1.

Following Theorem 1, there exist a family of $(\ell + M)-$dimensional vectors $\{\alpha_i : i \in S\}$ and a family of $(\ell+M)\times(\ell+M)$ matrices $\{\beta_{ij} : i,j \in S, \ i\ne j\}$ satisfying $\beta_{ij} = \beta_{ji}^{\mathrm{T}}$, such that

$$\theta_i(x^i) = \alpha_i + \sum_{j \ne i} \beta_{ij} B_j(x_j) \ , \ i \in S. \tag{5.2}$$

The families of potentials are given by

$$G_i(x_i) = \langle \alpha_i, B_i(x_i)\rangle + \log L'_i(x_i),$$
$$G_{ij}(x_i, x_j) = B_i^{\mathrm{T}}(x_i) \beta_{ij} B_j(x_j).$$

Let us note that two variables $X_i$ and $X_j$ are (spatially) conditionally independent, given $\{X_k,\ k \neq i,j\}$, if (and only if) $\beta_{ij} = 0$. In this case, we say that the sites $i$ and $j$ are neighbours. Thus the neighbourhood $\partial i$ of $i$ is the set $\{j\ :\ \beta_{ij} \neq 0\}$. Moreover we can substitute $x^i$ by $x_{\partial i}$ in the previous equations (5.1) and (5.2).

These auto-models for mixed state variables are completely and well defined once we choose admissible parameters $\{\alpha_i, \beta_{ij}\}$ that ensure the integrability condition (4.5).

### 5.2. Spatial cooperation behaviour

In many practical situations, we need to investigate the properties of local interactions of the system. Indeed, we want to know whether the field is spatially cooperative or spatially competitive (or neither of them). A standard definition of spatial cooperation (respectively, competition) is that at each site $i$, the conditional expectation $E\left[X_i \mid x^i\right]$ increases (respectively, decreases) with each neighbouring value $x_j$, $j \neq i$. For mixed-state auto-models valued in $E^S = (F \cup G)^S$, we must adapt these definitions. For each $i$ we define the function $x^i \longrightarrow R(x^i) = E\left[X_i 1_G(X_i) \mid x^i\right]$ and we study its variations in each coordinate $x_j$ of $x^i$, where $x_j \in G$ and $x_j \in \partial i$. Then we define spatial cooperation (or competition) similarly as the classical definition by substituting $R(x^i)$ for $E\left[X_i \mid x_{\partial i}\right]$. Let us note that this definition coincides with the classical one in the case $E = G$.

In the particular case where $E = \{0\} \cup (0, \infty)$ with $F = \{0\}$ and $G = (0, \infty)$, for any random mixed state variable $X$ on $E$ as defined in §2 with density function (2.2), we have

$$E[X] = (1 - \gamma)\int_0^\infty x\, g(x)\, dx = (1 - \gamma)E[X|1_{X>0}] = E[X 1_G(X)].$$

Then we conclude that for mixed state variables in $E^S = (\{0\} \cup (0, \infty))^S$, the generalised definition for spatial cooperation (competition) above meets the classical one. This will not occur anymore in the case of an atomic value different from zero.

### 5.3. A translation invariant and symmetric mixed-state auto-model with the four nearest neighbours system

Let us consider the four nearest neighbours system on a two-dimensional lattice, $S = [1, M] \times [1, N]$: each site $i \in S$ has four neighbours denoted by $\{i_e = i + (0, 1), i_w = i - (0, 1), i_n = i - (1, 0), i_s = i + (1, 0)\}$, with obvious neighbour adjustments near the boundary. We assume translation invariance in the sense that the parameters are functions of the displacement between sites; we assume spatial symmetry, which implies that the matrices $\beta_{ij} = \beta_{ji}$ are symmetric; we allow possible anisotropy between the horizontal and vertical directions. Under all these conditions and from the result above, there exist a $(\ell+M)$−dimensional vector $\alpha$ and two $(\ell+M) \times (\ell+M)$ symmetric matrices $\{\beta^{(1)}, \beta^{(2)}\}$ such that for all $i$, $\alpha_i = \alpha$, and, for all $\{i,j\}$, , $\beta_{ij} = 0$ unless $i$ and $j$ are neighbours, in which case

$$\beta_{i,i_e} = \beta_{i,i_w} = \beta^{(1)}, \quad \beta_{i,i_n} = \beta_{i,i_s} = \beta^{(2)}.$$

Moreover, the translation invariance implies that the local conditional density function $f_i$, hence the functionals $\theta_i$, $B_i$, $H'_i$ and $L'_i$ in (5.1), are independent of $i$. The potentials are given by

$$G_i(x_i) = G(x_i) = \langle \alpha, B(x_i) \rangle + \log L'(x_i),$$

$$G_{ij}(x_i, x_j) = G(x_i, x_j) = \begin{cases} B^{\mathrm{T}}(x_i)\beta^{(1)}B(x_j), & i-j = \pm(0,1), \\ B^{\mathrm{T}}(x_i)\beta^{(2)}B(x_j), & i-j = \pm(1,0), \\ 0, & \text{otherwise.} \end{cases}$$

The natural parameter of $f_i$ equals to

$$\theta_i(x^i) = \theta(x_{\partial i}) = \alpha + \beta^{(1)}\{B(x_{i_e}) + B(x_{i_w})\} + \beta^{(2)}\{B(x_{i_n}) + B(x_{i_s})\}. \quad (5.3)$$

### 5.4. Parameter estimation

It is well-known that the maximum likelihood method needs intensive computational approximations for Markov random fields. An efficient remedy relies on the pseudo-likelihood estimator introduced by [4]. Theoretical results for this estimator in the general framework of Markov random fields can be found in e.g. [8]. In the case of multi-parameter auto-models, [9] provide conditions under which this estimator is consistent. We refer the reader to this paper where this theory is developed in details. In particular, it applies to the present class of mixed-state auto-models. In the later §7, we will use this pseudo-likelihood estimator for the modelling of motion measurements from video sequences.

## 6. Mixed exponential auto-models

In this section, we focus on auto-models with mixed exponential conditional distributions. The relative simplicity of the model allows us a complete study of the various properties of the model, without getting bogged down in the parameters. Moreover, the exponential distribution itself is commonly used for modelling, e.g., records of pluviometry data.

From Besag's seminal paper [4], we know that several classical auto-models imply spatial competition between the neighbouring sites. This is particularly the case for the auto-exponential scheme. We will see that this fact appears again for mixed-state auto-model with exponential distributions. To overcome this limitation, we propose two alternatives by means of data truncation or data censoring.

### 6.1. Mixed state auto-models with exponential conditionals

We consider the mixed state space $E = \{0\} \cup (0, \infty)$, and we assume that the conditional distributions $f_i(x_i \mid \cdot)$ belong to the family of mixed state exponential distributions as defined in Example 1 of §2. We write

$$f_i(x_i \mid x^i) = f_i(x_i \mid x_{\partial i}) = \gamma_i(x_{\partial i})\delta(x_i) + (1 - \gamma_i(x_{\partial i}))\delta^*(x_i)g_{\lambda_i(x_{\partial i})}(x_i)$$
$$= H'_i(x_{\partial i}) \exp\langle \theta_i(x_{\partial i}), B(x_i)\rangle,$$

where the natural parameter and the sufficient statistics are, noting that $x\delta^*(x) = x$,

$$\theta_i(x_{\partial i}) = (\theta_{1,i}(x_{\partial i}), \theta_{2,i}(x_{\partial i}))^{\mathrm{T}} = \left(\log \frac{(1 - \gamma_i(x_{\partial i}))\lambda_i(x_{\partial i})}{\gamma_i(x_{\partial i})}, \lambda_i(x_{\partial i})\right)^{\mathrm{T}},$$
$$B(x) = (\delta^*(x), -x)^{\mathrm{T}}.$$

Here, the reference state is 0. Besides, obviously, the family of efficient statistics $B(x)$ verify Assumption 3. Therefore, following the previous general result for mixed state auto-models: there exist a family of vectors $\alpha_i = (a_i, b_i)^{\mathrm{T}}$ and $2 \times 2$ matrices $\beta_{ij} = \begin{pmatrix} c_{ij} & d_{ij} \\ f_{ij} & e_{ij} \end{pmatrix}$, with $\beta_{ij} = \beta_{ji}^{\mathrm{T}}$, such that the energy function equals

$$Q(x_1, \ldots, x_n) = \sum_{i \in S}(a_i\delta^*(x_i) - b_ix_i)$$
$$+ \sum_{\{i,j\}} \left(c_{ij}\delta^*(x_i)\delta^*(x_j) - d_{ij}\delta^*(x_i)x_j - f_{ij}x_i\delta^*(x_j) + e_{ij}x_ix_j\right).$$

We note that the model can be spatially asymmetrical if $d_{ij} \neq f_{ij}$. This is particularly interesting for our mixed state auto-models, where $\delta^*(x_i)x_j$ and $x_i\delta^*(x_j)$ may be interpreted as different situations. We can also think about models with oriented graphs.

The local natural parameters are

$$\theta_{1,i}(x_{\partial i}) = a_i + \sum_{\{i,j\}} (c_{ij}\delta^*(x_j) - d_{ij}x_j), \tag{6.1}$$

$$\theta_{2,i}(x_{\partial i}) = b_i + \sum_{\{i,j\}} (f_{ij}\delta^*(x_j) - e_{ij}x_j). \tag{6.2}$$

And we also have the reciprocal correspondence

$$\lambda_i(x_{\partial i}) = \theta_{2,i}(x_{\partial i}), \qquad \gamma_i(x_{\partial i}) = \frac{\theta_{2,i}(x_{\partial i})}{\theta_{2,i}(x_{\partial i}) + e^{\theta_{1,i}(x_{\partial i})}}.$$

It remains to make certain the well-definiteness of the model, that is to ensure the integrability condition (4.5). Necessarily, we must have for all $i$, $\gamma_i \in [0, 1]$ and $\lambda_i > 0$. Since the $x_j$'s belong to $[0, \infty)$, this leads to require the following conditions:

**(A)**  (i) for all $\{i,j\}$, $e_{ij} \leq 0$.
(ii) for all $i$ and any subset $A \subset S\setminus\{i\}$, $b_i + \sum_{j\in A} f_{ij} > 0$ (in particular $b_i > 0$).

Fortunately enough, these necessary conditions also ensure the integrability condition (4.5).

**Proposition 1.** *Under Conditions (A), the auto-model with mixed exponential conditionals is well-defined.*

*Proof.* The configuration space $\Omega$ can be decomposed as

$$\Omega = \sum_{A\subset S} \Omega_A, \quad \text{with} \quad \Omega_A = \{x:\ x_i > 0,\ i \in A;\quad x_i = 0,\ i \notin A\}.$$

We have

$$\int_\Omega \exp Q(x)\nu(dx) = \sum_{A\subset S} \int_{\Omega_A} \exp Q(x)\ m^{\otimes S}(dx),$$

with $m(dx) = \delta_0(dx) + \lambda(dx)$. Therefore, Condition (4.5) holds if and only if

$$\forall A \subset S,\ \int_{\Omega_A} \exp Q(x)\ m^{\otimes S}(dx) < +\infty.$$

Moreover for $x \in \Omega_A$, with $f_{ij} = d_{ji}$,

$$Q(x) = \sum_{i\in A}[a_i - b_i x_i] + \sum_{\{i,j\}\subset A}(1,-x_i)\beta_{ij}(1,-x_j)^{\mathrm{T}}$$

$$= \sum_{i\in A} a_i + \sum_{i,j\in A,\langle i,j\rangle} c_{ij} - \sum_{i\in A} x_i\left[(b_i + \sum_{j\in A:\langle i,j\rangle} f_{ij}) - \sum_{j\in A:\langle i,j\rangle} e_{ij}x_j\right].$$

As $e_{ij} \leq 0$, we have for some constant $C > 0$, and still $x \in \Omega_A$,

$$C^{-1}Q(x) \leq -\sum_{i\in A}\left(b_i + \sum_{j\in A,\ \langle i,j\rangle} f_{ij}\right)x_i.$$

Let $|A| = \mathrm{card}(A)$. By Conditions **(A)**, $b_i + \sum_{j\in A,\ \langle i,j\rangle} f_{ij} > 0$ and we finally obtain

$$\int_{\Omega_A} \exp Q(x)\ m^{\otimes S}(dx) = \int_{(0,\infty)^{|A|}} \exp Q(x) \prod_{i\in A} \lambda(dx_i) < \infty.$$

The proof is complete. □

Let us notice that in the context of $n$ "ordinary" variables, [3] claim that Condition (A) is both necessary and sufficient for $n \geq 2$. This is true for $n = 2$ but the condition is not necessary anymore if $n \geq 3$.

Let us examine local interactions between neighbouring sites. Considering the generalised definition of spatial cooperation (respectively, competition) for mixed-state auto-models given in §5.2 , we are looking for the variations of

$$E\left[X_i 1_{X_i>0} \mid x_{\partial i}\right] = (1 - \gamma_i(x_{\partial i})) \frac{1}{\theta_{2,i}(x_{\partial i})}$$
$$= \frac{1}{1 + \theta_{2,i}(x_{\partial i}) \exp\{-\theta_{1,i}(x_{\partial i})\}} \frac{1}{\theta_{2,i}(x_{\partial i})}.$$

Under Conditions (**A**), in particular $e_{ij} \leq 0$, we see that the parameter $\theta_{2,i}(x_{\partial i})$ defined in (6.2) is an increasing function of neighbouring values $x_j > 0$. As $E[X_i 1_{X_i>0} \mid x_{\partial i}]$ is a decreasing function of $\theta_{2,i}(x_{\partial i})$, we conclude that the model cannot be spatially cooperative, although the precise dependence of the other parameter $\theta_{1,i}(x^i)$ in the $x_j$'s will vary according to the values of the $\{d_{ij}\}$'s. Similarly to the auto-exponential scheme of [4], this locally non cooperative behaviour seems inappropriate in many application fields.

To overcome this drawback, there are two commonly used approaches, namely data truncation and data censoring. We adapt below these two methods for mixed exponential auto-models.

### 6.2. Cooperative mixed exponential auto-models by truncation

First let us define a mixed truncated exponential variable $X$. The state space is $E = \{0\} \cup (0, K]$ where $K$ is a given (arbitrary) positive constant. The continuous component on $(0, K]$ of $X$ follows a truncated exponential distribution with the probability density function

$$g_\lambda(x) = H(\lambda) e^{-\lambda x} \mathbf{1}_{(0,K]}(x), \quad H(\lambda) = \frac{\lambda}{1 - e^{-\lambda K}}.$$

Thus, the probability density function of $X$ equals

$$f_\theta(x) = H'(\theta) \exp\langle \theta, B(x) \rangle$$

with

$$\theta = (\theta_1, \theta_2)^{\mathrm{T}} = \left(\log \frac{(1-\gamma)H(\lambda)}{\gamma}, \lambda\right)^{\mathrm{T}}, \quad B(x) = (\delta^*(x), -x \mathbf{1}_{(0,K]}(x))^{\mathrm{T}}.$$

Note that conversely we have

$$\lambda = \theta_2, \quad \gamma = \frac{\theta_2}{\theta_2 + e^{\theta_1}(1 - e^{-\theta_2 K})}.$$

Let us consider a mixed-state auto-model for $X = \{X_i, i \in S\}$, whose conditional distributions lay in the family of the mixed truncated exponential distributions above. Here the reference state is 0 and the family of sufficient statistics $B$ verify Assumption 3. By Theorem 1, there exist a family of vectors $\alpha_i = (a_i, b_i)^{\mathrm{T}}$ and $2 \times 2$ matrices $\beta_{ij} = \begin{pmatrix} c_{ij} & d_{ij} \\ f_{ij} & e_{ij} \end{pmatrix}$, with $\beta_{ij} = \beta_{ji}^{\mathrm{T}}$, such that the energy function equals

$$Q(x_1,\ldots,x_n) = \sum_{i \in S}(a_i,b_i)B(x_i) + \sum_{\{i,j\}} B(x_i)^{\mathrm{T}} \beta_{ij} B(x_j).$$

Because of the truncation, $\exp Q$ is always integrable.

The natural parameters of local conditional distributions are written

$$\theta_i(x_{\partial i}) = \begin{pmatrix} \theta_{i,1}(x_{\partial i}) \\ \theta_{i,2}(x_{\partial i}) \end{pmatrix} = \begin{pmatrix} a_i + \sum_{j \in \partial i} \left\{ c_{ij}\delta^*(x_j) - d_{ij}x_j 1_{(0,K]}(x_j) \right\} \\ b_i + \sum_{j \in \partial i} \left\{ f_{ij}\delta^*(x_j) - e_{ij}x_j 1_{(0,K]}(x_j) \right\} \end{pmatrix} \quad (6.3)$$

As for the conditions on the parameters, we keep the requirement: $\theta_{i,2}(x_{\partial i}) > 0$ (which implies $\gamma_i(x_{\partial i}) \in (0,1)$). This is clearly satisfied under the following

*Assumption* 4. For all $i \in S$, $\quad b_i + \sum_{j \in \partial i} \min(0, f_{ij}, f_{ij} - e_{ij}K) > 0$.

To understand whether the system is spatially cooperative or not, let us first examine, for the mixed truncated exponential variable $X$, the variation of $E[X1_{(0,K]}(X)]$ with respect to its parameters $\{\theta_1, \theta_2\}$. If we denote $Z$ a random variable following a truncated exponential distribution with the density $g_\lambda$, a simple calculus leads to $E[Z] = K\left(\frac{1}{\lambda K} - \frac{1}{e^{\lambda K}-1}\right)$, which decreases from $\frac{1}{2}K$ to 0 when $\lambda$ raises from 0 to $\infty$. On the other side, $1-\gamma$ is decreasing with respect to $\theta_2$ and increases with $\theta_1$. Finally, $E[X1_{(0,K]}(X)] = (1-\gamma)E(Z)$ is decreasing in $\theta_2$ and increasing in $\theta_1$.

Gathering this result together with (6.3), we deduce the variation of $E[X_i 1_{(0,K]}(X_i) \mid x_{\partial i}]$, relatively to the neighbouring values $x_j$, $j \in \partial i$ which are positives. Let us introduce the following assumptions.

*Assumption* 5. For all $i,j \in S$, $d_{ij} \leq 0$, $e_{ij} \geq 0$.

*Assumption* 6. For all $i,j \in S$, $d_{ij} \geq 0$, $e_{ij} \leq 0$.

We thus have proved the following

**Proposition 2.** *Assume Assumption* 4 *holds. Then,*

(i) *The mixed auto-model with mixed truncated exponential conditionals is well-defined.*
(ii) *The model is spatially cooperative under Assumption* 5.
(iii) *The model is spatially competitive under Assumption* 6.

Therefore, this family of auto-models can exhibit spatial cooperation as well as spatial competition.

Let us give an application to the translation invariant and symmetric scheme with the four nearest neighbours system, as introduced in §5.3, and with the additional assumption of isotropy. Then the parameters are $\alpha = (a,b)^T$ and $\beta^{(1)} = \beta^{(2)} = \begin{pmatrix} c & d \\ d & e \end{pmatrix}$.

In this case, Assumptions 4 and 5 reduce to the conditions

$$d \leq 0, \quad e \geq 0, b + 4(d - eK) > 0,$$

which makes the model spatially cooperative. Similarly, the model is spatially competitive under Assumptions 4 and 6 which bring down to

$$d \geq 0, \quad e \leq 0, b > 0.$$

Note that in the case of spatial cooperation (Assumptions 4 and 5), the positive parameters $e_{ij}$'s in $\theta_{i,2}(x_{\partial i})$, Eq. (6.3), give a measure of the strength of the spatial cooperation: the bigger are the values of these parameters, the stronger is the spatial cooperation realized in the model. However, Assumptions 4 and 5 imply $0 \leq e_{ij} < h_{ij}/K$ for some positive constants $h_{ij}$. Therefore if the truncation level $K$ is large, the implied spatial cooperation becomes limited.

### 6.3. Cooperative mixed exponential auto-models by censoring

As previously, let $K$ be a fix positive constant. We consider the mixed censored exponential variable $X$ defined in §3.2. Let us recall the expression of the corresponding probability density function.

$$f_\theta(x) = H'(\theta) \exp\langle \theta, B(x) \rangle$$

with

$$\begin{aligned}
\theta &= (\theta_1, \theta_2, \theta_3)^{\mathrm{T}} = \left( \log \frac{\alpha}{(1-\alpha)} + \lambda K, \log \lambda + \lambda K, \lambda \right)^{\mathrm{T}}, \\
B(x) &= (1_{\{0\}}(x), 1_{(0,K)}(x), -x 1_{(0,K)}(x))^{\mathrm{T}}.
\end{aligned}$$

The reference state is $K$, and we notice that the components of $\theta$ are dependent. Conversely we have

$$\lambda = \theta_3, \qquad \alpha = \frac{e^{\theta_1 - \theta_3 K}}{1 + e^{\theta_1 - \theta_3 K}}.$$

Let us consider now a mixed-state auto-model for $X = \{X_i, \ i \in S\}$, whose conditional distributions belong to the family of mixed censored exponential distributions above. For sake of simplicity, we assume spatial symmetry. Applying again Theorem 1, there exist a family of 3-dimensional vectors $\alpha_i = (r_i, a_i, b_i)^{\mathrm{T}}$ and $3 \times 3$ symmetric matrices

$$\beta_{ij} = \begin{pmatrix} s_{ij} & u_{ij} & t_{ij} \\ u_{ij} & c_{ij} & d_{ij} \\ t_{ij} & d_{ij} & e_{ij} \end{pmatrix}.$$

such that the energy function equals to

$$Q(x_1, \ldots, x_n) = \sum_{i \in S} (r_i, a_i, b_i) B(x_i) + + \sum_{\{i,j\}} B(x_i)^{\mathrm{T}} \beta_{ij} B(x_j).$$

The natural parameters of local conditional distributions are

$$\theta_i(x_{\partial i}) = \begin{pmatrix} \theta_{i,1}(x_{\partial i}) \\ \theta_{i,2}(x_{\partial i}) \\ \theta_{i,3}(x_{\partial i}) \end{pmatrix} = \begin{pmatrix} r_i \\ a_i \\ b_i \end{pmatrix} + \sum_{j \in \partial i} \begin{pmatrix} s_{ij} & u_{ij} & t_{ij} \\ u_{ij} & c_{ij} & d_{ij} \\ t_{ij} & d_{ij} & e_{ij} \end{pmatrix} \begin{pmatrix} 1_{\{0\}}(x_j) \\ 1_{(0,K)}(x_j) \\ -x_j 1_{(0,K)}(x_j) \end{pmatrix}$$
(6.4)

Because of the censoring, $\exp Q$ is always integrable, but the parameters have to meet the requirement $\theta_{i,3}(x_{\partial i}) > 0$. This is clearly satisfied under the assumption

*Assumption 7.* For all $i \in S$, $\quad b_i + \sum_{j \in \partial i} \min(t_{ij},\ d_{ij},\ d_{ij} - e_{ij}K) > 0$.

We turn now to the study of the possible cooperative or competitive features of the model. We first examine, for the mixed censored exponential variable $X$, the variation of $E[X 1_{(0,K]}(X)]$ with respect to its parameters $\{\theta_1, \theta_2, \theta_3\}$. If we denote $Z'$ a random variable following a censored exponential distribution, we have
$$E[Z'] = \frac{1 - e^{-\lambda K}}{\lambda},$$
which decreases from $K$ to $0$ when $\lambda$ raises from $0$ to $\infty$. Therefore,
$$E[X 1_{(0,K]}(X)] = (1-\alpha) E(Z') = \frac{1}{1 + e^{\theta_1 - \theta_3 K}} \frac{1 - e^{-\theta_3 K}}{\theta_3}.$$

If the parameters $\theta_1$ and $\theta_3$ depend on some real variable $\eta$ in such a way that both functions $\theta_3$ and $\theta_1 - \theta_3 K$ are decreasing (respectively, increasing) in $\eta$, then the expectation $E[X 1_{(0,K]}(X)]$ is an increasing (respectively, decreasing) function in $\eta$.

Coming back to $E[X_i 1_{(0,K]}(X_i) \mid x_{\partial i}]$, from (6.4), we have
$$\theta_{i,3}(x_{\partial i}) = b_i + \sum_{j \in \partial i} \{t_{ij} 1_{\{0\}}(x_j) + d_{ij} 1_{(0,K)}(x_j) - e_{ij} x_i 1_{(0,K)}(x_j)\},$$

and
$$\begin{aligned}(\theta_{i,1} - \theta_{i,3} K)(x_{\partial i}) =& (r_i - b_i K) \\ &+ \sum_{j \in \partial i} \{(s_{ij} - t_{ij} K) 1_{\{0\}}(x_j) + (u_{ij} - d_{ij} K) 1_{(0,K)}(x_j) \\ &- (t_{ij} - e_{ij} K) x_i 1_{(0,K)}(x_j)\}.\end{aligned}$$

Let us introduce the following assumptions.

*Assumption 8.* For all $i, j \in S$, $e_{ij} \geq 0$, $\quad t_{ij} - e_{ij} K \geq 0$.

*Assumption 9.* For all $i, j \in S$, $e_{ij} \leq 0$, $\quad t_{ij} - e_{ij} K \leq 0$.

**Proposition 3.** *Assume Assumption 7 holds.*

(i) *The mixed auto-model whose local conditional distributions are mixed censored exponential distributions is well-defined.*
(ii) *The model is spatially cooperative under Assumption 8.*
(iii) *The model is spatially competitive under Assumption 9.*

The proof is straightforward, directly derived from the definition of spatial cooperation given in 5.2 and the above expressions of the parameters.

This family of auto-models is able to produce spatial cooperation as well as spatial competition, and this under rather light conditions, which makes the model attractive. To convince, let us look at the translation invariant and symmetric scheme with the four nearest neighbours system of §5.3, assuming also spatial isotropy. Therefore, the unique interaction matrix equals

$$\beta^{(1)} = \beta^{(2)} = \begin{pmatrix} s & u & t \\ u & c & d \\ t & d & e \end{pmatrix}.$$

In this case, Assumptions 7 and 8 reduce to the conditions

$$e \geq 0, \quad t - eK \geq 0, \quad b + 4\min(t, d - eK) > 0,$$

that makes the model spatially cooperative. Analogously, Assumptions 7 and 9 reduce to the conditions

$$e \leq 0, \quad t - eK \leq 0, \quad b + 4\min(t, d) > 0,$$

which guarantee the spatial competition behaviour of the model.

## 7. An application to motion analysis from video image sequences

### 7.1. Motion measurements from video sequences

Motion computation and analysis are of central importance in image analysis. Let $\{I_i(t)\}$ be an image sequence where $i = (i_1, i_2) \in S$ denotes the pixel locations and $t = 1, \ldots, T$ the time instants in the sequence. Roughly speaking, a motion map at time $t$, $X(t) = \{X_i(t)\} = \{\|v_i(t)\|\}$ is defined as the norm of the underlying motion field $\{v_i(t)\}$ which is estimated by a "regularised" minimisation of the sum of squares $\sum_i [I_{i+v_i(t)}(t+1) - I_i(t)]^2$. Usually some local smoothing procedures are needed to get a more robust motion map and we refer to [7] for details of these computations.

We consider here video sequences of natural scenes. Figure 1 displays three sample images from each of two sequences involving a moving escalator and trees under wind respectively. The corresponding motion maps $\{X_i(t)\}$ are displayed in Figure 2. Next, sample histograms from these motion maps are presented in Figure 3. As a matter of fact, these histograms present a composite picture. An important peak appears at the origin accounting for regions where no motion is present, while a continuous component encompasses actual motion magnitudes in the images.

### 7.2. An mixed-state auto-model with positive Gaussian distributions

We follow the general construction of mixed states auto-models of §4. First, we consider a positive mixed-state Gaussian variable $X$, defined in Example 3 of §2. Then $X$ has the following density function :

$$f_\theta(x) = \gamma\delta(x) + (1-\gamma)\delta^*(x)g_\xi(x) = \exp\left[\langle\theta, B(x)\rangle + \log\gamma\right], \qquad (7.1)$$

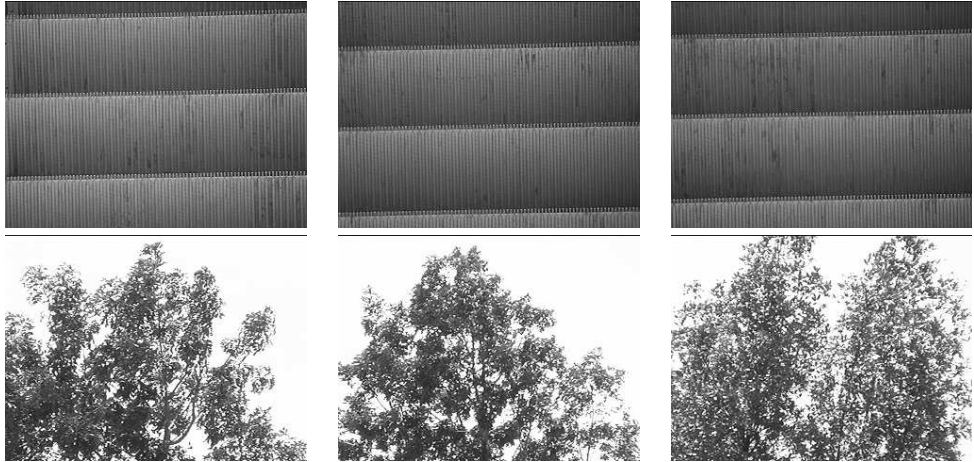

FIGURE 1. *Sample images from two videos. Top row: a moving escalator; bottom row: trees.*

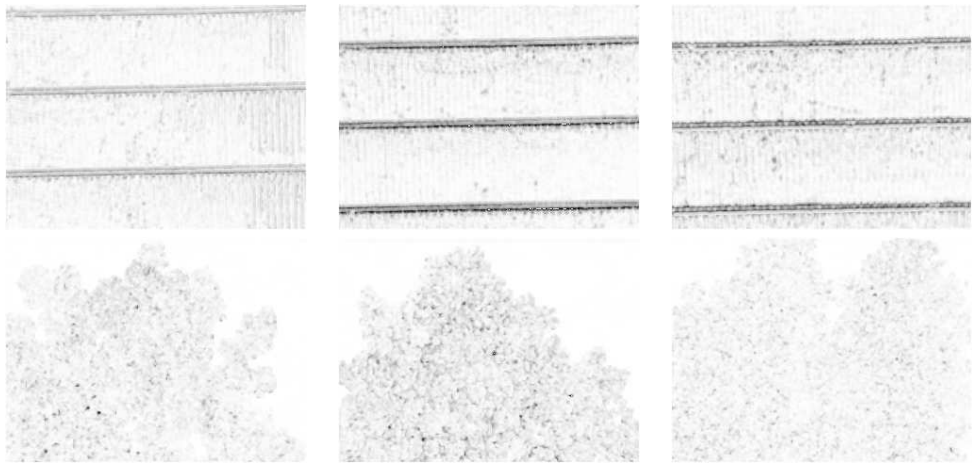

FIGURE 2. *Sample motion measures $\{X_i(t)\}$ from the videos of Figure 1. Top row: a moving escalator; bottom row: a tree (white=0; black=maximum value).*

where $\xi = (2\sigma^2)^{-1}$, $g_\xi(x) = 2(2\pi\sigma^2)^{-1/2} \exp\{-\frac{1}{2\sigma^2}x^2\}$, and

$$\theta = (\theta_1, \theta_2)^{\mathrm{T}} = \left(\log \frac{(1-\gamma)g_\xi(0)}{\gamma}, \xi\right)^{\mathrm{T}}, \quad B(x) = (\delta^*(x), -x^2)^{\mathrm{T}}. \qquad (7.2)$$

To construct auto-models for the motion maps observations $\{X_i(t)\}$, we assume that the family of conditional distributions $f_i(x_i|x^i)$ belongs to the family of mixed positive Gaussian distribution given in (7.1). By Theorem 1, there exist a family of vectors $\alpha_i = (a_i, b_i)^{\mathrm{T}} \in \mathbb{R}^2$ and a family of $2 \times 2$ matrices $\{\beta_{ij}\}$

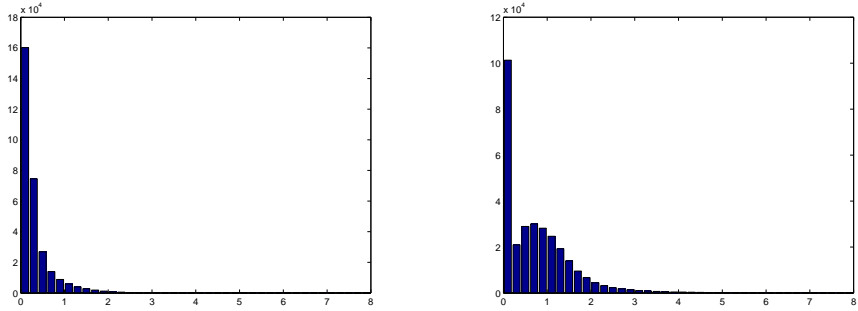

Figure 3. *Sample histograms of motion measures $\{X_i(t)\}$ of Figure 2.*

satisfying $\beta_{ij} = \beta_{ji}^{\mathrm{T}}$, such that

$$\theta_i(x^{(i)}) = \alpha_i + \sum_{j \neq i} \beta_{ij} B(x_j) . \qquad (7.3)$$

Moreover the associated energy function is given by

$$Q(x_1, \ldots, x_n) = \sum_{i \in S} \left[ a_i \delta^*(x_i) - b_i x_i^2 \right] + \sum_{\{i,j\}} (\delta^*(x_i), -x_i^2) \beta_{ij} (\delta^*(x_j), -x_j^2)^{\mathrm{T}} . \qquad (7.4)$$

To analyse the motion measurements, we consider the specification of §5.3, namely a translation invariant and spatially symmetric auto-model with the four-nearest-neighbours system and possible anisotropy between the horizontal and vertical directions. Then the parametrisation reduces to one vector $\alpha = (a, b)^{\mathrm{T}}$ and two $2 \times 2$ matrices $\beta^{(1)}$ and $\beta^{(2)}$ such that $\forall i,\ \alpha_i = \alpha$, and $\forall \{i, j\},\ \beta_{ij} = 0$ unless $i$ and $j$ are neighbours in which case

$$\beta_{i,i_e} = \beta^{(1)} = \begin{pmatrix} c_1 & d_1 \\ d_1 & e_1 \end{pmatrix} = \beta_{i_w,i}^{\mathrm{T}}, \quad \beta_{i,i_n} = \beta^{(2)} = \begin{pmatrix} c_2 & d_2 \\ d_2 & e_2 \end{pmatrix} = \beta_{i_s,i}^{\mathrm{T}}.$$

Moreover, the present context asks for spatial cooperation and we need further to constrain the parameters $d_k$ and $e_k$, $k = 1, 2$ to be zero. The resulting automodel has four parameters $\phi = (a, b, c_1, c_2)$ and is well-defined (admissible) under an unique condition: $b > 0$. Here, we use the pseudo-likelihood method to estimate these parameters.

Let us mention that in the context of image segmentation, Salzenstein and Pieczynski [12] have previously proposed a fuzzy image segmentation model where the fuzzy labels are a particular instance of mixed-state variables with values in $[0, 1]$.

### 7.3. *Experiments*

The experiments are conducted in order to evaluate whether the model above can correctly account for a fundamental characteristic of an homogeneous tex-ture, namely spatial isotropy or spatial anisotropy. For the present four-nearest-neighbours Gaussian mixed-state model, the spatial isotropy occurs if (and only if) $c_1 = c_2$.

We fit this model to several motions maps like those displayed in Figure 2. First we consider motions from trees (bottom row of the figure). A typical set of parameter estimates is $\hat{\phi} = (\hat{a}, \hat{b}, \hat{c}_1, \hat{c}_2) = (-5.805, 3.044, 3.057, 2.954)$. The parameters $c_1$ and $c_2$ are almost identical with regard to standard deviations of these estimates computed at other time instants of the same tree sequence. Therefore, the believed spatial isotropy for these motions is well reflected here.

Next we consider the motion maps from a moving escalator (top row of Figure 2). Since the motion is a vertical one, we clearly have anisotropy. A typical set of parameter estimates is $\hat{\phi} = (\hat{a}, \hat{b}, \hat{c}_1, \hat{c}_2) = (-6.512, 0.320, 2.192, 3.598)$. Therefore, the difference between $c_1$ and $c_2$ appears to be significant and the mixed-state model is able to reflect the spatial anisotropy of the considered motion. More experiments on motion analysis can be found in [5].

## Acknowledgement

This work originated from several discussions with Xavier Guyon. The authors are also grateful to Gwënaelle Piriou for the experiments provided in Section 7.